\documentclass[leqno,12pt]{article}
\usepackage{amsmath,amssymb,color,a4wide}
\usepackage{MnSymbol}
\usepackage{tikz-cd}
\usepackage{xr-hyper}
\usepackage{hyperref}

\usepackage[curve,matrix,arrow,cmtip]{xy}
\usepackage{verbatim}
\NoComputerModernTips

\newdir^{ (}{{}*!/-3pt/\dir^{(}}
\newdir_{ (}{{}*!/-3pt/\dir_{(}}

\def\bigtimes{\mathop{\raisebox{-2pt}{\huge $\times$\kern-1pt}}}

\def\OM{\mathchoice
  {\rlap{\kern3.2pt$\overline{\phantom{L}}$}M}
  {\rlap{\kern3.2pt$\overline{\phantom{L}}$}M}
  {\rlap{\kern2.4pt$\scriptstyle\overline{\phantom{L}}$}M}
  {\rlap{\kern1.8pt$\scriptscriptstyle\overline{\phantom{L}}$}M}}

\let\le\leqslant
\let\ge\geqslant
\let\leq\leqslant
\let\geq\geqslant

 \def\Cl{\mathop{\rm Cl}\nolimits}

 \def\deg{\mathop{\rm deg}\nolimits}

\def\ord{\mathop{\rm ord}\nolimits}

\def\GL{\mathop{\rm GL}\nolimits}
\def\SL{\mathop{\rm SL}\nolimits}

\def\Sp{\mathop{\rm Sp}\nolimits}

\def\Id{{\rm Id}}

\let\phi\varphi
\let\theta\vartheta
\let\epsilon\varepsilon
\let\setminus\smallsetminus

\newtheorem{Thm}{Theorem}[section]
\newtheorem{Prop}[Thm]{Proposition}
\newtheorem{Lem}[Thm]{Lemma}

\newtheorem{Cor}[Thm]{Corollary}

\newtheorem{Def}[Thm]{Definition}

\newtheorem{Rem}[Thm]{Remark}
\newtheorem{Ex}[Thm]{Example}

\def\UseTheoremCounterForNextEquation{\setcounter{equation}{\value{Thm}}\addtocounter{Thm}{1}}

\def\qed{{\hskip0pt\unskip\unskip\nobreak\hfil\penalty50
          \hskip1em\hbox{}\nobreak\hfil
           {$\square$}
          \parfillskip=0pt\finalhyphendemerits=0
          \par}\medskip}

\newenvironment{Proof}
               {\noindent{\bf Proof.}\ }
               {\qed}
							
               {\noindent{\bf Proof Sketch.}\ }
               {\qed}

               {\noindent{\bf Proof of #1.}\ }
               {\qed}

\newcommand{\BC}{{\mathbb{C}}}

\newcommand{\BF}{{\mathbb{F}}}

\newcommand{\BP}{{\mathbb{P}}}
\newcommand{\BQ}{{\mathbb{Q}}}

\newcommand{\BZ}{{\mathbb{Z}}}

\newcommand{\CE}{{\cal E}}

\newcommand{\CM}{{\cal M}}
\newcommand{\CN}{{\cal N}}

\newcommand{\CS}{{\cal S}}
\newcommand{\CT}{{\cal T}}

\newcommand{\CW}{{\cal W}}

\newbox\mybox
\def\arrover#1{\mathrel{
       \setbox\mybox=\hbox spread 1.4em
              {\hfil$\scriptstyle#1$\hfil}
       \vbox{\offinterlineskip\copy\mybox
             \hbox to\wd\mybox{\rightarrowfill}}}}

\def\larrover#1{\mathrel{
       \setbox\mybox=\hbox spread 1.4em
              {\hfil$\scriptstyle#1\vphantom{g}$\hfil}
       \vbox{\offinterlineskip\copy\mybox
             \hbox to\wd\mybox{\leftarrowfill}}}}

\def\ontoover#1{\mathrel{
       \setbox\mybox=\hbox spread 1.4em
              {\hfil$\scriptstyle#1\vphantom{g}$\hfil}
       \vbox{\offinterlineskip\copy\mybox
             \hbox to\wd\mybox{\rightarrowfill\hskip-2.8mm
                               $\rightarrow$}}}}
\def\leftontoover#1{\mathrel{
       \setbox\mybox=\hbox spread 1.4em
              {\hfil$\scriptstyle#1\vphantom{g}$\hfil}
       \vbox{\offinterlineskip\copy\mybox
             \hbox to\wd\mybox{$\leftarrow$\hskip-2.8mm
                               \leftarrowfill}}}}
\let\longto\longrightarrow

\def\Cinf{{\BC}_\infty}
\def\Finf{F_\infty}

\def\Bigskip{\bigskip\bigskip}

\newbox\dotDdbox
\newbox\dotDtbox
\newbox\dotDsbox
\newbox\dotDssbox
\setbox\dotDdbox=\vbox{\hbox{\kern3.5pt\bf.}\vskip-4.5pt\hbox{$D$}}
\setbox\dotDtbox=\vbox{\hbox{\kern3.5pt\bf.}\vskip-4.5pt\hbox{$D$}}
\setbox\dotDsbox=\vbox{\hbox{\kern2.1pt.}
                       \vskip-6.6pt\hbox{$\scriptstyle D$}}
\setbox\dotDssbox=\vbox{\hbox{\kern1.6pt.}
                        \vskip-8.1pt\hbox{$\scriptscriptstyle D$}}

\newcommand{\OmegaModGammaU}{\Gamma_U\backslash \Omega^r}

\newcommand{\Spm}{\text{Spm}\;}

\newcommand{\abs}[1]{\left\lvert #1\right\rvert}

\newcommand{\Proposition}{Proposition }

\newcommand{\Definition}{Definition }

\newcommand{\PropositionCite}{Prop.\,}

\newcommand{\TheoremCite}{Thm.\,}

\newcommand{\norm}{\xi}   

\begin{document}

\title{Drinfeld modular forms of arbitrary rank\\
Part I: Analytic Theory}

\author{Dirk Basson \and Florian Breuer$^{1,2}$ \and Richard Pink$^2$}

\footnotetext[1]{Supported by the Alexander von Humboldt foundation, and by the NRF grant BS2008100900027.}
\footnotetext[2]{Supported through the program ``Research in Pairs'' by Mathematisches Forschungsinstitut Oberwolfach in 2010.}

\date{May 27, 2018}
\maketitle

\centerline{To Lukas}
\bigskip

\begin{abstract}
This is the first of a series of articles providing a foundation for the theory of Drinfeld modular forms of arbitrary rank~$r$. In the present part, we develop the analytic theory. Most of the work goes into defining and studying the $u$-expansion of a weak Drinfeld modular form, whose coefficients are weak Drinfeld modular forms of rank $r-1$. Based on that we give a precise definition of when a weak Drinfeld modular form is holomorphic at infinity and thus a Drinfeld modular form in the proper sense.
\end{abstract}

{\advance\baselineskip by -6pt
\tableofcontents
}


\externaldocument{Part_2}
\externaldocument{Part_3}
\setcounter{section}{0}

\Bigskip


\newpage
\noindent {\Large\bf Introduction}
\bigskip

\noindent
In \cite{Drinfeld1}, Drinfeld introduced elliptic modules, now called \emph{Drinfeld modules}, in order to prove a special case of the Langlands conjectures for $\GL_2$ over function fields. These objects share many properties with elliptic curves, though their rank can be an arbitrary integer $r\ge1$. In particular, Drinfeld constructed a moduli space of Drinfeld modules of rank $r$ with a suitable level structure, both as an algebraic variety and with an analytic uniformisation as a quotient of an $r-1$ dimensional symmetric space~$\Omega^r$. This $\Omega^r$ is a rigid analytic space over a field $\Cinf$ of positive characteric and plays the role of the complex upper half plane. In the case $r=2$ Drinfeld \cite{Drinfeld2} used automorphic forms on $\Omega^r$ with values in $\BQ_\ell$ to prove a case of the Langlands conjectures for the associated automorphic representations on $\GL_2$. 

But there is also a natural definition of \emph{Drinfeld modular forms} on $\Omega^r$ with values in the field $\Cinf$ of positive characteristic. Goss \cite{GossES} was the first to explicitly refer to these, defining them both algebraically, in the way Katz did in \cite{Katz}, and analytically as (rigid analytic) holomorphic functions on~$\Omega^r$. In the case $r=2$, where these are functions of one variable, it was relatively straightforward to write down the necessary condition of \emph{holomorphy at infinity}. This led to the development of a theory of Drinfeld modular forms of rank~$2$, for instance by Gekeler \cite{GekelerDMC};
see \cite{gekeler_survey_1999} for a survey.

For $r\ge3$ the situation concerning holomorphy at infinity is more subtle. In the related case of Siegel modular forms of genus $\ge2$ the problem disappears, because the necessary condition at infinity holds automatically by the K\"ocher principle. One explanation for this is the fact that in the Satake compactification of the Siegel moduli space of abelian varieties the boundary has codimension $\ge2$. By contrast, the moduli space of Drinfeld modules is always affine, so in any compactification as an algebraic variety the boundary has codimension~$1$; hence a condition at infinity is always required.

That condition is important for several reasons. On the one hand many relevant modular forms that one can construct naturally, such as Eisenstein series, satisfy it automatically. On the other hand a condition at infinity is necessary for one of the main structural results, the fact that the space of modular forms of given level and weight is finite dimensional.

The condition at infinity can be expressed in two quite different ways. The analytic way says that the $u$-expansion (which is a kind of Fourier expansion) of a modular form consists only of terms with non-negative index. For the algebraic way one identifies the analytic modular forms with sections of an invertible sheaf on a moduli space. Then one requires a compactification of this moduli space as a projective algebraic variety over~$\Cinf$ together with an extension of the invertible sheaf. The crucial step is to prove that a modular form satisfies the analytic condition at infinity if and only if the corresponding section on the moduli space extends to a section on that compactification. The finite dimensionality is then a direct consequence of the fact that the space of sections of a coherent sheaf on a projective algebraic variety is always finite dimensional.
Using the \emph{Satake compactification} of a Drinfeld moduli space, the third author \cite{Pink} has already established much of the necessary algebro-geometric theory for this.


\medskip
The present paper is the first of a series of articles together with \cite{BBP2} and \cite{BBP3}, whose aim is to provide the rest of the theory and thereby a foundation for the theory of Drinfeld modular forms of arbitrary rank. The present Part I develops the basic analytic theory, including $u$-expansions and holomorphy at infinity. Part II will identify the analytic modular forms discussed here with the algebraic modular forms defined in \cite{Pink} and deduce qualitative consequences such as the finite dimensionality of the space of modular forms of given level and weight. Part III will illustrate the general theory by constructing and studying some important families of modular forms.

\medskip

For a discussion on the history of Drinfeld modular forms of higher rank, see \cite[\S7]{BassonBreuer2017}. We mention here also two recent papers by Gekeler \cite{GekelerHigherRankI,GekelerHigherRankII} in which he constructs the building map from $\Omega^r$ to the Bruhat-Tits building of $\GL_r$ and uses this to study the growth and vanishing behaviours of important families of modular forms for $\GL_r(\BF_q[t])$.

\subsubsection*{Outline  of this paper}

In Section \ref{Sec:Notation} we introduce our notation and define the Drinfeld period domain $\Omega^r$ with its action by $\GL_r(F)$ for a global function field $F$. Weak modular forms for an arithmetic subgroup $\Gamma < \GL_r(F)$ are defined as holomorphic functions from $\Omega^r$ to $\Cinf$ satisfying the functional equation (\ref{DefOfActionOnF}) linking $f(\gamma(\omega))$ to $f(\omega)$ for every $\gamma\in\Gamma$.

Further preparations are made in the next two sections. In Section \ref{Sec:Exp} we collect basic properties of exponential functions associated to strongly discrete subgroups of $\Cinf$, and we outline the rigid analytic structure of $\Omega^r$ in Section \ref{Sec:RigidStructure}.

Based on our choice of coordinates on $\Omega^r$, we identify a \emph{standard boundary component}, whose translates by $\GL_r(F)$ form all boundary components of codimension~$1$. Thus a weak modular form is holomorphic at all boundary components if and only if all its translates by $\GL_r(F)$ are holomorphic at the standard boundary component. The holomorphy at the standard boundary component is tested using the expansion with respect to a certain parameter $u$.

%

This parameter is defined in Section \ref{Sec:HolomorphyA}: We decompose elements $\omega\in\Omega^r$ as $\omega=\binom{\omega_1}{\omega'}$, where $\omega_1\in\Cinf$ is the first coordinate of $\omega$ and $\omega'\in\Omega^{r-1}$ consists of the remaining coordinates. Next, we assign to $\Gamma$ a subgroup $\Lambda'\subset F^{r-1}$ isomorphic to the subgroup $\Gamma_U < \Gamma$ of translations which fix $\omega'$. Then $\Lambda'\omega'\subset\Cinf$ is a strongly discrete subgroup and we can form the associated exponential function $e_{\Lambda'\omega'}$. Now $e_{\Lambda'\omega'}(\omega_1)$ is invariant under the translations $\Gamma_U$ and we define our parameter as its reciprocal $u:=u_{\omega'}(\omega_1)=e_{\Lambda'\omega'}(\omega_1)^{-1}$ in (\ref{DefOfuParameter}).

In Definition \ref{def:NbhdOfInfty} we define neighbourhoods of infinity in $\Omega^r$, then Theorem \ref{Thm:Parameter} states that the map $\binom{\omega_1}{\omega'}\mapsto \binom{u}{\omega'}$ induces rigid analytic isomorphisms from quotients of neighbourhoods of infinity by $\Gamma_U$ to so-called pierced tubular neighbourhoods in $\Cinf^\times\times\Omega^{r-1}$. 

This allows us to show in Section \ref{Sec:HolomorphyB} that every weak modular form $f$ admits a $u$-expansion 
\[
f(\omega) = \sum_{n\in\BZ}f_n(\omega')u_{\omega'}(\omega_1)^n
\]
converging on a neighbourhood of infinity (Proposition \ref{ModFormsLaurentExpansion}), whose coefficients $f_n$ are themselves weak modular forms on $\Omega^{r-1}$ (Theorem \ref{Thm:CoefficientsAreModularForms}). These are the main results of this paper.

Finally, we define modular forms in Section \ref{Sec:AMF} as weak modular forms all of whose rotations by elements of $\GL_r(F)$ admit $u$-expansions with terms of non-negative index. It follows from Propositions \ref{Prop:OrderInvariantB} and \ref{Prop:cosets} that this condition only needs to be tested for finitely many elements of $\GL_r(F)$. It will be shown in Part II of this series that this definition agrees with the algebraic definition of modular forms in \cite{Pink}.

\section{Weak modular forms}
\label{Sec:Notation}

Throughout this paper we fix a global function field $F$ of characteristic $p>0$, with exact field of constants $\BF_q$ of cardinality~$q$. We fix a place $\infty$ of $F$ and let $A$ denote the ring of elements of $F$ which are regular away from~$\infty$. 
This is a Dedekind domain with finite class group and group of units $A^\times=\BF_q^\times$.
Let $\pi\in F$ be a uniformising parameter at $\infty$, so that $|\pi|=q^{-\deg\infty}$. Let $\Finf\cong\BF_{q^{\deg\infty}}((\pi))$  denote the completion of $F$ at $\infty$, and $\Cinf$ the completion of an algebraic closure of $\Finf$.

We fix an unspecified non-zero constant $\norm \in \Cinf^\times$, whose value can be set for normalisation purposes. For example, if $F=\BF_q(t)$ and $A=\BF_q[t]$, there are certain advantages in letting $\norm$ be a period of the Carlitz module. For more general function fields $F$, a natural choice of $\norm$ is a period of a certain sign-normalised rank-one Drinfeld module, see \cite[Chapter IV (2.14) and (5.1)]{GekelerDMC}. However, we will not explicitly need the uniformisation in this article, so the reader loses nothing by assuming that $\norm = 1$.

The \emph{Drinfeld period domain of rank $r\ge 1$ over $F_\infty$} is usually defined as the set of points ${(\omega_1:\ldots:\omega_r)} \in \BP^{r-1}(\Cinf)$ for which $\omega_1,\ldots,\omega_r$ are $\Finf$-linearly independent. Any such point possesses a unique representative with $\omega_r=\norm$. We shall only work with these representatives, so we identify $\Omega^r$ with the following subset of $\Cinf^r$:
\UseTheoremCounterForNextEquation
\begin{equation}\label{DefOfOmega}
\Omega^r := \bigl\{ (\omega_1,\ldots,\omega_r)^T \in \Cinf^r \bigm|
\omega_1,\ldots,\omega_r\ \Finf\hbox{-linearly independent and\ }\omega_r=\norm
\bigr\}.
\end{equation}
We write the elements of $\Omega^r$ as $r \times 1$ matrices, i.e. column vectors.

For any point $\omega\in\Omega^r$ and any matrix $\gamma\in\GL_r(\Finf)$, the matrix product $\gamma\omega$ is again a column vector with $\Finf$-linearly independent entries. In particular its last entry is non-zero. Defining
\UseTheoremCounterForNextEquation
\begin{equation}\label{DefOfj}
j(\gamma,\omega) := \norm^{-1}\cdot(\hbox{last entry of\ }\gamma\omega) \in \Cinf^\times,
\end{equation}
we thus find that
\UseTheoremCounterForNextEquation
\begin{equation}\label{DefOfActionOnOmega}
\gamma(\omega) := j(\gamma,\omega)^{-1}\gamma\omega
\end{equation}
again lies in~$\Omega^r$. This defines a left action of $\GL_r(\Finf)$ on~$\Omega^r$. Also, for any $\gamma$, $\delta\in\GL_r(\Finf)$ a direct calculation shows that
\UseTheoremCounterForNextEquation
\begin{equation}\label{jCocycle}
j(\gamma\delta,\omega) = j\big(\gamma,\delta(\omega)\big) j(\delta,\omega).
\end{equation}

For any function $f : \Omega^r \to \Cinf$ and any integers $k$ and $m$ we define the function $f|_{k,m}\gamma : \Omega^r\to\Cinf$ by
\UseTheoremCounterForNextEquation
\begin{equation}\label{DefOfActionOnF}
(f|_{k,m}\gamma)(\omega)
:= \det(\gamma)^{m}j(\gamma,\omega)^{-k}f\big(\gamma(\omega)\big).
\end{equation}
By direct calculation we deduce from (\ref{jCocycle}) that
\UseTheoremCounterForNextEquation
\begin{equation}\label{fCocycle}
(f|_{k,m}\gamma\delta)(\omega) = ((f|_{k,m}\gamma)|_{k,m}\delta)(\omega).
\end{equation}
Thus (\ref{DefOfActionOnF}) defines a right action of $\GL_r(\Finf)$ on the space of all functions $f:\Omega^r \to \Cinf$.

For later use note also that, if $\gamma=a\cdot\Id_r$ for the identity matrix $\Id_r\in\GL_r(F)$, then $j(\gamma,\omega)=a$ and $\gamma(\omega)=\omega$ and $\det(\gamma)=a^r$; and hence 
\UseTheoremCounterForNextEquation
\begin{equation}\label{IdrAction}
f|_{k,m}(a\cdot\Id_r) = a^{rm-k}f.
\end{equation}

\begin{Rem}\label{ConventionsRemA}
\rm There are different conventions about whether $\Omega^r$ consists of row or column vectors and about how $\GL_r(\Finf)$ acts on it. For instance, the first and third authors \cite{basson_thesis}, \cite{Pink} like Drinfeld \cite{Drinfeld1} use row vectors and the action $(\gamma,\omega) \mapsto \omega\gamma^{-1}$, whereas Kapranov \cite{Kapranov} uses column vectors and the action by left multiplication $(\gamma,\omega) \mapsto \gamma\omega$. These conventions differ not only by transposition, but also by the outer automorphism $\gamma\mapsto(\gamma^T)^{-1}$ of $\GL_r$.
The present article uses column vectors and left multiplication in order to make things compatible with the existing literature on rank $2$ Drinfeld modular forms.
\end{Rem}

The set $\Omega^r$ can be endowed with the structure of a rigid analytic space. Experts may
be content with the fact that $\Omega^r$ is an admissible open subset of $\BP^{r-1}(\Cinf)$ 
and inherits its rigid analytic structure, while others may consult Section \ref{Sec:RigidStructure} 
for more details. A holomorphic function on $\Omega^r$ is a global section of the structure 
sheaf of $\Omega^r$, but a more useful characterisation  is that a function $f: \Omega^r\to\Cinf$ is holomorphic
if and only if it is a uniform limit of rational functions on $\BP^{r-1}(\Cinf)$ whose poles 
all lie outside $\Omega^r$. 

\begin{Def}\label{Def:WeakModForm}
Consider integers $k$ and $m$ and an arithmetic subgroup $\Gamma<\GL_r(F)$. A \emph{weak 
modular form of weight $k$ and type $m$ for $\Gamma$} is a holomorphic function $f:\Omega^r\to\Cinf$ 
which for all $\gamma\in\Gamma$ satisfies
$$f|_{k,m}\gamma = f.$$
The space of these functions is denoted by $\CW_{k,m}(\Gamma)$.
\end{Def}

Since $\Gamma$ is an arithmetic subgroup of $\GL_r(F)$, its determinant $\det(\Gamma)$ is 
a finite subgroup of $F^\times$ and therefore contained in the multiplicative group of 
the field of constants~$\BF_q^\times$. Thus its order is a divisor of $q-1$, and the 
definition depends only on $m$ modulo this divisor; in other words we have
\UseTheoremCounterForNextEquation
\begin{equation}\label{TypeCongruence}
\CW_{k,m}(\Gamma) = \CW_{k,m'}(\Gamma)
\hbox{\ whenever $m\equiv m'$ modulo $|\det(\Gamma)|$.}
\end{equation}
On the other hand, for any $\alpha\in\BF_q^\times$ we have $f|_{k,m}(\alpha\cdot\Id_r) = \alpha^{rm-k}f$ by (\ref{IdrAction}). Thus 
\UseTheoremCounterForNextEquation
\begin{equation}\label{TypeZero}
\CW_{k,m}(\Gamma) = 0
\hbox{\ unless $k \equiv rm$ modulo $\bigl|\Gamma\cap\{\hbox{scalars}\}\bigr|$.}
\end{equation}

In the case $m=0$ we will suppress all mention of $m$ and abbreviate $f|_k\gamma := f|_{k,m}\gamma$
 and $\CW_k(\Gamma) := \CW_{k,m}(\Gamma)$. By (\ref{TypeCongruence}) we may always do this when 
$\Gamma<\SL_r(F)$.

\medskip
For later use we note the following direct consequence of (\ref{fCocycle}):

\begin{Prop}\label{Prop:WConj}
For any $\delta\in\GL_r(F)$ we have $f\in \CW_{k,m}(\Gamma)$ if and only if $f|_{k,m}\delta \in \CW_{k,m}(\delta^{-1}\Gamma\delta)$.
\end{Prop}

In general the space $\CW_{k,m}(\Gamma)$ is infinite dimensional. A finite dimensional subspace 
of `non-weak' modular forms will be characterised by conditions at infinity. The formulation of 
these conditions requires some preparation in the next two sections.

\section{Exponential functions}
\label{Sec:Exp}

A subgroup $H\subset\Cinf$ is called {\em strongly discrete} if its intersection with every ball of finite radius is finite. For any such subgroup, any $z\in\Cinf$, and any $\epsilon>0$, there are at most finitely many elements $h\in H\setminus\{0\}$ with $\abs{\frac{z}{h}}\ge\epsilon$. In this case the product
\UseTheoremCounterForNextEquation
\begin{equation}\label{ExpDef}
e_H(z) := z\cdot\!\!\prod_{h\in H\smallsetminus\{0\}}\left(1-\frac{z}{h}\right)
\end{equation}
converges in~$\Cinf$, defining the {\em exponential function $e_H : \Cinf \to \Cinf$ associated to~$H$}.

\begin{Prop}\label{exp1}
For any strongly discrete subgroup $H\subset\Cinf$, the function $e_H : \Cinf \to \Cinf$ is 
holomorphic, surjective, and has simple zeros at the points in~$H$ and no other zeros. It induces an 
isomorphism of additive groups and rigid analytic spaces
$$\Cinf/H\stackrel{\sim}{\longrightarrow}\Cinf.$$
The function $e_H$ possesses an everywhere convergent power series expansion
$$e_H(z) = \sum_{i=0}^\infty e_{H,p^i}z^{p^i}$$
with $e_{H,p^i}\in\Cinf$ and $e_{H,1}=1$. If $H$ is an $\BF_q$-subspace, the expansion has the form
$$e_H(z) = \sum_{j=0}^\infty e_{H,q^j}z^{q^j}.$$
If $H$ is finite, then $e_H(z)$ is a polynomial of degree $|H|$ in~$z$.
\end{Prop}

\begin{Proof}
When $H\subset\Cinf$ is an $A$-lattice (see below), this is proved in \cite[\S4.2]{GossBS} and 
\cite[Prop. 1.27]{GossES}. The case where $H\subset\Cinf$ is merely a strongly discrete subgroup follows in exactly the same way.
%
%
\end{Proof}


%

\begin{Prop}\label{exp2}
\begin{enumerate}
\item[(a)] For any two strongly discrete subgroups $H'\subset H\subset\Cinf$, the subgroup $e_{H'}(H)\subset\Cinf$ is strongly discrete and isomorphic to $H/H'$, and we have
$$e_H = e_{e_{H'}(H)}\circ e_{H'}.$$
\item[(b)] For any strongly discrete subgroup $H\subset\Cinf$ and any $a\in\Cinf^\times$, the subgroup $aH\subset\Cinf$ is strongly discrete, and we have
$$e_{aH}(az) = a e_H(z).$$
\end{enumerate}
\end{Prop}

\begin{Proof}
For (a) see \cite[(1.12)]{GekelerPowerSums}, and (b) follows immediately from the definition.
\end{Proof}

An {\em $A$-lattice of rank $r$ in $\Cinf$} is a strongly discrete projective $A$-submodule $\Lambda\subset\Cinf$ of rank $r$.


\begin{Prop}\label{exp3}
Let $H\subset\Cinf$ be an $A$-lattice of rank~$r$. Then for any $a\in A$ there exists a 
unique $\BF_q$-linear polynomial $\phi^H_a(z)$ of degree $|a|^r$ satisfying
\[ \phi^H_a(e_H(z)) = e_H(az)\]
for all $z\in\Cinf$. The map $\phi^H: a\mapsto \phi^H_a$ is a Drinfeld $A$-module of rank~$r$.
\end{Prop}
\begin{Proof}
\cite[\TheoremCite 4.3.1]{GossBS}
\end{Proof}

\section{The rigid analytic structure of $\Omega^r$}
\label{Sec:RigidStructure}

Throughout the following we denote by $B(0,\rho) := \{z\in\Cinf : |z|\le \rho\}$ the closed disk of radius $\rho>0$ centred at~$0$, and by $B(0,\rho)'=B(0,\rho)\setminus\{0\}$ the associated punctured disk. We will also consider the annulus $D(0,\sigma,\rho) := {\{z\in\Cinf\mid \sigma\le|z|\le \rho\}}$. Note that $B(0,\rho)$ and $D(0,\sigma,\rho)$ are affinoid whenever $\sigma,\rho\in|\Cinf^\times|$.

\medskip
We will describe the rigid analytic structure of $\Omega^r$ by covering it by suitable affinoid subspaces. Two such coverings already appear in \cite{Drinfeld1}, and one of them is described in more detail in \cite{SchneiderStuhler}. We follow the approach in \cite{SchneiderStuhler}, but adapt it to our convention that $\omega_r=\norm$.

We say that a linear form $F_\infty^r\to F_\infty$ is \emph{unimodular} if its largest coefficient has absolute value 1. 
For any $\Finf$-rational hyperplane $H\subset\BP^{r-1}(\Cinf)$, we choose a unimodular linear form $\ell_H$ that defines it.  Then $|\ell_H(\omega)|$ is well-defined and non-zero for any $\omega\in\Omega^r$. 
Using the standard norm $|\omega|:=\max_{1\le i\le r} |\omega_i|$ on~$\Cinf^r$, we set 
\UseTheoremCounterForNextEquation
\begin{equation}\label{DOmega}
h(\omega)\ :=\ \tfrac{1}{|\omega|}\cdot\inf\{|\ell_H(\omega)|:H \text{ an }\Finf\text{ hyperplane}\},
\end{equation}
which measures the distance from $\omega\in\Omega^r$ to all boundary components combined. For any $n\in\BZ^{>0}$ we also define
\UseTheoremCounterForNextEquation
\begin{equation}\label{Omega_n}
\Omega^r_n:=\{\omega\in\Omega^r: h(\omega)\ge |\pi|^n\}.
\end{equation}
Since $|\pi| < 1$, these subsets satisfy $\Omega^r_1 \subset \Omega^r_2 \subset\ldots$ and their union is~$\Omega^r$.

\begin{Lem}\label{Lem:Omeganbounds}
Every $\omega\in\Omega_n^r$ satisfies $|\xi|\leq |\omega| \leq |\xi||\pi|^{-n}$.
\end{Lem}

\begin{Proof}
The first inequality follows from $\omega_r = \xi$. 
Next, since $\omega\mapsto\omega_r$ is a unimodular $\Finf$-linear form, (\ref{DOmega}) implies that $|\omega|h(\omega) \leq |\xi|$, from which the second inequality follows.
\end{Proof}

\begin{Prop}\label{AffinoidCovering}
For each $n\in\BZ^{>0}$, the set $\Omega^r_n$ is an affinoid subdomain of $\BP^{r-1}(\Cinf)$. Together they form an admissible covering of $\Omega^r$, endowing it with the structure of an admissible open subset of $\BP^{r-1}(\Cinf)$.
\end{Prop}

\begin{Proof} See version (C) of the proof of \cite[\PropositionCite 1]{SchneiderStuhler}.\end{Proof}

Using the second (finer) covering in \cite[\PropositionCite 6.2]{Drinfeld1}, Drinfeld showed that any arithmetic subgroup $\Gamma<\GL_r(F)$ acts \emph{discontinuously} on $\Omega^r$, which means that there exists an admissible covering $\Omega^r=\bigcup_{i\in I}U_i$ such that for each $i\in I$, the set $\{\gamma\in\Gamma\,|\,\gamma(U_i)\cap U_i\ne\emptyset\}$ is finite. 
Thus, for every subgroup $G<\Gamma$, the quotient $G\backslash \Omega^r$ exists as a rigid analytic space (see \cite[\S6.4]{FresnelvdPut}).

\medskip
For the following recall that a function $f:U\to\Cinf$ on an admissible open subset $U\subset \Omega^r$ is called \emph{holomorphic} if it is a section of the sheaf of functions on this space, or equivalently, if it is a uniform limit $f=\lim_{n\to\infty}f_n$ of rational functions $f_n:\BP^{r-1}(\Cinf)\dashedrightarrow \Cinf$ with no poles in~$U$.

\medskip
In the next section we shall need bounds on the values of certain exponential functions 
when we restrict to $\omega\in\Omega^r_n$. For this we require the following estimates:

\begin{Lem}\label{Lem:BoundingAbsoluteValues}
For any $\gamma\in\GL_r(\Finf)$ there exist positive constants $c_1$, $c_2$ and $c_3$ 
such that for every $\omega\in\Omega^r$ we have
\begin{enumerate}
\item[(a)] $h(\omega) \leq c_1 |j(\gamma,\omega)||\omega|^{-1}\le 1$;
\item[(b)] $\displaystyle |\gamma(\omega)|\le c_2 h(\omega)^{-1}$; and
\item[(c)] $\displaystyle h(\gamma(\omega)) \ge c_3 h(\omega)$. 
\end{enumerate}
\end{Lem}

\begin{Proof}
Let $x$ be an entry of the last row of $\gamma$ of maximal absolute value, 
and set $c_1 := |x^{-1}\norm|>0$. Then by the definition (\ref{DefOfj}) of $j(\gamma,\omega)$, the value $x^{-1}\norm j(\gamma,\omega)$ is a unimodular $\Finf$-linear combination of the $\omega_i$'s, so we obtain
\[h(\omega) \leq \frac{|x^{-1}\norm j(\gamma,\omega)|}{|\omega|} \leq 1.\]
This proves (a).

Next, let $c_2'$ be the largest absolute value of an entry of $\gamma$. Then the matrix product satisfies $|\gamma\omega| \leq c_2'|\omega|$ and so 
$|\gamma(\omega)| = |j(\gamma,\omega)^{-1}||\gamma\omega| \leq |j(\gamma,\omega)^{-1}|c_2'|\omega| \leq c_1c_2'h(\omega)^{-1}$, where the last inequality follows from (a). This proves (b) with $c_2 = c_1c_2'$.

For (c) let $c_3'$ denote the largest absolute value of an entry of $\gamma^{-1}$. Let $\ell$ be an arbitrary unimodular $\Finf$-linear form, which we write as a row vector, so that $\ell(\omega)=\ell\omega$. Choose $m_\ell\in\Cinf^\times$ such that $\ell_\gamma:=m_\ell\ell\gamma$ is a unimodular linear form. Then the entries of $m_\ell\ell = m_\ell\ell \gamma\cdot\gamma^{-1}$ have absolute value at most~$c_3'$; hence $|m_\ell| \leq c_3'$. Since $\gamma(\omega)=j(\gamma,\omega)^{-1}\gamma\omega$, using the linearity of $\ell$ and the definition of $h(\omega)$ we find that
$$\frac{|\ell(\gamma(\omega))|}{|\gamma(\omega)|}
\ =\ \frac{|\ell\gamma\omega|}{|\gamma\omega|}
\ =\ \frac{|m_\ell|^{-1}|\ell_\gamma\omega|}{|\gamma\omega|}
\ \ge\ \frac{c_3^{\prime-1}|\ell_\gamma\omega|}{c_2'|\omega|}
\ \ge\ \frac{h(\omega)}{c_2'c_3'}.$$
Varying $\ell$, the definition of $h(\gamma(\omega))$ now implies (c) with $c_3 := (c_2'c_3')^{-1}$.
\end{Proof}

\section{Neighbourhoods of infinity}
\label{Sec:HolomorphyA}

From now on we assume that $r\ge2$. Let $U$ denote the algebraic subgroup of $\GL_{r,F}$ 
of matrices of the form
\UseTheoremCounterForNextEquation
\begin{equation}\label{DefU}
\left(\begin{array}{c|c}
1 & *\ \ldots\ * \\ \hline
\begin{array}{c} 0 \\[-5pt] \vdots \\[-4pt] 0 \end{array}
& \Id_{r-1}
\end{array}\right)
\end{equation}
where $\Id_{r-1}$ denotes the identity matrix of size $(r-1)\times(r-1)$.
Fix an arithmetic subgroup $\Gamma<\GL_r(F)$ and set
\UseTheoremCounterForNextEquation
\begin{equation}\label{DefGammaU}
\Gamma_U := \Gamma\cap U(F).
\end{equation}
Then for all $\gamma\in\Gamma_U$ and $\omega\in\Omega^r$ we have $\det(\gamma)=j(\gamma,\omega)=1$; 
hence every weak modular form for $\Gamma$ is a $\Gamma_U$-invariant function on~$\Omega^r$.

Viewing elements of $F^{r-1}$ as $1\times(r-1)$-matrices (row vectors), consider the isomorphism
\UseTheoremCounterForNextEquation
\begin{equation}\label{eq:iota}
\iota: F^{r-1} \stackrel{\sim}{\longto} U(F),\ \ v'\mapsto
\left(\begin{matrix}1 & v' \\ 0 & \Id_{r-1} \end{matrix}\right).
\end{equation}
Since $\Gamma$ is commensurable with $\GL_r(A)$, the subgroup
\UseTheoremCounterForNextEquation
\begin{equation}\label{DefLambdaPrime}
\Lambda' := \iota^{-1}(\Gamma_U) \subset F^{r-1}
\end{equation}
is commensurable with $A^{r-1}$. On the other hand, recall that $\Omega^r$ is the set of column vectors $\omega = (\omega_1,\ldots,\omega_r)^T\in\Cinf^r$ with $F_\infty$-linearly independent entries and $\omega_r=\norm$. For any such $\omega$ we have  $\omega' := (\omega_2, \ldots,\omega_r)^T \allowbreak \in \Omega^{r-1}$, hence
$$\Omega^r \subset \Cinf\times \Omega^{r-1}$$
inside $\Cinf^r = \Cinf\times\Cinf^{r-1}$. Accordingly we write $\omega = \binom{\omega_1}{\omega'}$. 
The definition (\ref{DOmega}) then directly implies that $h(\omega)\le h(\omega')$ and hence $\Omega^r_n \subset\Cinf\times\Omega^{r-1}_n$.

For any element $\lambda'\in\Lambda'$ we can form the matrix product $\lambda'\omega' \in \Cinf$. The definition (\ref{DefOfActionOnOmega}) of the action on $\Omega^r$ then implies that
\UseTheoremCounterForNextEquation
\begin{equation}\label{Arrival}
\textstyle \iota(\lambda') (\binom{\omega_1}{\omega'}) = \binom{\omega_1 + \lambda'\omega'}{\omega'},
\end{equation}
which extends the action to $\Cinf\times\Omega^{r-1}$ by the same formula.
For any $\omega'\in\Omega^{r-1}$ observe that $\Lambda'\omega' := \{\lambda'\omega'\mid \lambda'\in\Lambda'\}$ is a strongly discrete subgroup of $\Cinf$, because $\Lambda'$ is commensurable with $A^{r-1}$ and the entries of $\omega'$ are $\Finf$-linearly independent. Thus the function
\UseTheoremCounterForNextEquation
\begin{equation}\label{eq:TerraNova}
\textstyle \Cinf\times \Omega^{r-1}\to \Cinf,\ 
\left[\binom{\omega_1}{\omega'}\right] \mapsto e_{\omega'\Lambda'}(\omega_1)
\end{equation}
is well-defined and $\Gamma_U$-invariant. 

As usual in a metric space, for any point $z\in\Cinf$ and any subset $X\subset\Cinf$ we write $d(z,X) := \inf\{|z-x|:x\in X\}$. Then we have:

\begin{Prop}\label{NewEarth}
\begin{enumerate}
\item[(a)] The function (\ref{eq:TerraNova}) is holomorphic. 
\item[(b)] For any $n\in\BZ^{>0}$ and $c>0$ there exists a constant $c_n>0$, such that for any $\omega'\in\Omega^{r-1}_n$ and any $\omega_1\in\Cinf$ with $|\omega_1|<c$ we have $|e_{\Lambda'\omega'}(\omega_1)|<c_n$.
\item[(c)] For any $n\in\BZ^{>0}$ and $R_n>0$ there exists a constant $c_n>0$, such that for any $\omega'\in\Omega^{r-1}_n$ and any $\omega_1\in\Cinf$ with $d(\omega_1,F_\infty^{r-1}\omega')<R_n$ we have $|e_{\Lambda'\omega'}(\omega_1)|<c_n$.
\item[(d)] For any $\omega'\in\Omega^{r-1}$ and $\omega_1\in\Cinf$ we have $|e_{\Lambda'\omega'}(\omega_1)|\ge d(\omega_1,F_\infty^{r-1}\omega')$.
\end{enumerate}
\end{Prop}

\begin{Proof}
The function is defined by the product $e_{\omega'\Lambda'}(\omega_1) = \omega_1\cdot\prod_{0\not=\lambda'\in\Lambda'}(1-\frac{\omega_1}{\lambda'\omega'})$, whose factors we examine in turn.
First, as $\Lambda'\subset F^{r-1}$ is commensurable with $A^{r-1}$, there exists a constant $a\in A\setminus\{0\}$ with $\Lambda'\subset a^{-1}A^{r-1}$. 
Next observe that any $\lambda'\in\Lambda'\setminus\{0\}$ determines a unimodular $F_\infty$-linear form $\frac{\lambda'}{|\lambda'|}$. For any $\omega'\in\Omega_n^{r-1}$ it follows that
\UseTheoremCounterForNextEquation
\begin{equation}\label{eq:intermediate}
|\lambda'\omega'|
\ =\ |\lambda'|\cdot\bigl|\tfrac{\lambda'}{|\lambda'|}\omega'\bigr|
\ \smash{\stackrel{\eqref{DOmega}}{\geq}}\ 
|\lambda'|\cdot h(\omega')\cdot|\omega'|
\ \smash{\stackrel{\eqref{Omega_n}}{\geq}}\ 
|\lambda'|\cdot|\pi^n|\cdot|\omega'|
\ \smash{\stackrel{\ref{Lem:Omeganbounds}}{\geq}}\ 
|\lambda'|\cdot|\pi^n|\cdot|\xi|.
\end{equation}
As $\lambda'$ runs through $\Lambda'\setminus\{0\}$, the value $|\lambda'\omega'|$ thus goes to $\infty$ uniformly over~$\Omega_n^{r-1}$. Varying $n$ this implies that the function is holomorphic, proving (a).

To prove (b) observe next that all factors $1-\frac{\omega_1}{\lambda'\omega'}$ with $|\lambda'\omega'|\ge|\omega_1|$ have absolute value less than or equal to~$1$. Since now $|\omega_1|< c$, we deduce that
\UseTheoremCounterForNextEquation
\begin{equation}\label{eq:product}
|e_{\omega'\Lambda'}(\omega_1)|\ <\ c\cdot\!\!\prod_{0<|\lambda'\omega'|<c}\frac{c}{|\lambda'\omega'|}.
\end{equation}
Since $\Lambda'\subset a^{-1}A^{r-1}$, for any $\lambda'\in\Lambda'\setminus\{0\}$ we have $|\lambda'|\ge\frac{1}{|a|}$. From (\ref{eq:intermediate}) we thus deduce that $|\lambda'\omega'| \geq \frac{|\pi^n\xi|}{|a|}$. In particular each factor in the product (\ref{eq:product}) satisfies $\frac{c}{|\lambda'\omega'|}\leq \frac{c|a|}{|\pi^n\xi|}$; hence it is bounded by a constant that is independent of~$\omega'$. 
On the other hand, if $|\lambda'\omega'|<c$, the inequality (\ref{eq:intermediate}) implies that $|\lambda'|<\frac{c}{|\pi^n\xi|}$. Thus each coefficient of $a\lambda'\in A^{r-1}$ has absolute value $<\frac{c|a|}{|\pi^n\xi|}$, the number of possibilities for which is bounded independently of~$\omega'$. The number of factors in (\ref{eq:product}) is thus also bounded independently of~$\omega'$, and so is therefore the total value of the product, proving (b).

To prove (c) write $\omega_1 = x\omega'+y$, where $x\in\Finf^{r-1}$ and $y\in\Cinf$ with $|y|<R_n$. Since $\Lambda'\subset F^{r-1}$ is commensurable with $A^{r-1}$, the factor group $\Finf^{r-1}/\Lambda'$ is compact. Thus there exists a constant $\alpha>0$ depending only on $A$ and~$\Lambda'$, such that every $x\in\Finf^{r-1}$ can be written in the form $x=\lambda'+x_0$ for some $\lambda'\in\Lambda'$ and $x_0\in\Cinf$ with $|x_0|<\alpha$. Together we then have $\omega_1 = \lambda'\omega'+(x_0\omega'+y)$ with $|x_0\omega'|<\alpha|\omega'| \leq \alpha|\xi\pi^{-n}|$ by Lemma \ref{Lem:Omeganbounds} and hence $|x_0\omega'+y| < \max\{\alpha|\xi\pi^{-n}|,R_n\}$. By part (b) this implies that $|e_{\Lambda'\omega'}(\omega_1)| = |e_{\Lambda'\omega'}(x_0\omega'+y)| < c_n$ for some constant $c_n>0$ that is independent of $\omega_1$ and $\omega'$, proving (c).

To prove (d) write $\omega_1 = \lambda'_0\omega'+y$ with $\lambda'_0\in\Lambda'$ and $y\in\Cinf$ such that $|y|$ is minimal. Then for all $\lambda'\in\Lambda'$ we have $|y-\lambda'\omega'|\ge|y|$. 
If $|y|\ge|\lambda'\omega'|$, this implies that $|y-\lambda'\omega'|\ge|\lambda'\omega'|$ and hence $\bigl|1-\frac{y}{\lambda'\omega'}\bigr|\ge1$.
If $|y|<|\lambda'\omega'|$, we directly deduce that $\bigl|1-\frac{y}{\lambda'\omega'}\bigr|=1$.
Writing $e_{\omega'\Lambda'}(\omega_1) = e_{\omega'\Lambda'}(y) = y\prod_{0\not=\lambda'\in\Lambda'}(1-\frac{y}{\lambda'\omega'})$, we conclude that all factors in the product satisfy $\bigl|1-\frac{y}{\lambda'\omega'}\bigr|\ge1$. Thus it follows that $|e_{\omega'\Lambda'}(\omega_1)| \ge |y| \ge d(\omega_1,F_\infty^{r-1}\omega')$, proving (d).
\end{Proof}

\begin{Prop}\label{Prop:OmegaByGammaUIso}
The action of $\Gamma_U$ on $\Cinf\times\Omega^{r-1}$ from (\ref{Arrival}) is free and discontinuous and the quotient $\Gamma_U\backslash(\Cinf\times\Omega^{r-1})$ exists as a rigid analytic space. Moreover we have an isomorphism of rigid analytic spaces
$$\textstyle\CE:\ 
\Gamma_U\backslash(\Cinf\times \Omega^{r-1})\longto \Cinf\times\Omega^{r-1},\ \ 
\left[\binom{\omega_1}{\omega'}\right] \mapsto 
\binom{e_{\omega'\Lambda'}(\omega_1)}{\omega'}.$$ 
\end{Prop}

\begin{Proof}
The subsets $U_{\rho,n} := B(0,\rho)\times\Omega^{r-1}_n$ for all $\rho\in|\Cinf^\times|$ and $n>0$ form an admissible affinoid covering of $\Cinf\times\Omega^{r-1}$. For fixed $\rho$ and~$n$, consider a non-trivial element $\gamma=\iota(\lambda')\in\Gamma_U$ such that $\gamma U_{\rho,n}\cap U_{\rho,n} \neq \emptyset$. Choose $\binom{\omega_1}{\omega'}\in U_{\rho,n}$ such that $\gamma\binom{\omega_1}{\omega'} = \binom{\omega_1+\lambda'\omega'}{\omega'}\in U_{\rho,n}$. Then $|\lambda'\omega'| \le \max\{|\omega_1|,|\omega_1+\lambda'\omega'|\} \leq\rho$. Since $\frac{\lambda'}{|\lambda'|}$ is a unimodular $F_\infty$-linear form, we obtain
$$\frac{\rho}{|\lambda'|}
\ \geq\ \frac{|\lambda'\omega'|}{|\lambda'|} 
\ \stackrel{\eqref{DOmega}}{\geq}\ h(\omega')|\omega'| 
\ \stackrel{\eqref{Omega_n}}{\geq}\ |\pi|^n|\omega'| 
\ \stackrel{\ref{Lem:Omeganbounds}}{\geq}\ |\pi|^n|\xi|.$$
Thus we have $|\lambda'| \leq \rho|\pi|^{-n}|\xi|^{-1}$, where the right hand side depends only on $\rho$ and~$n$. Since $\lambda'$ lies in~$\Lambda'$, which is commensurable with $A^{r-1}$, this leaves only finitely many possibilities for $\gamma = \iota(\lambda')$. Thus the action is discontinuous. The action is clearly free, and the existence of the quotient as rigid analytic space follows from \cite[\S6.4]{FresnelvdPut}.

By Proposition \ref{exp1} we obtain a well-defined bijective and holomorphic map~$\CE$. As the derivative of $e_{\Lambda'\omega'}(X)$ with respect to $X$ is identically 1, the morphism is also \'etale. By \Proposition \ref{Prop:EtaleIso} below it is therefore an isomorphism.
\end{Proof}

\begin{Prop}
\label{Prop:EtaleIso}
Let $f:X\to Y$ be a morphism of rigid analytic spaces defined over an algebraically closed field $K$ which is \'etale and bijective. Then $f$ is an isomorphism.
\end{Prop}

\begin{Proof}
(The proof is based on the analogous argument for schemes at \cite[Tag 02LC]{stacks-project}.)
First we show that $f$ is universally injective, i.e., for any morphism $g:Y'\to Y$ the morphism $f':X' := X\times_Y Y'\to Y'$ is injective. So consider any points $x', x''\in X'$ mapping to the same point $y'\in Y'$. Then they also map to the same point $y\in Y$, and by the bijectivity of $f$ they therefore also map to the same point $x\in X$. Thus $x'$ and $x''$ lie in the fiber product $x\times_yy'$ which, since all these points have the same residue field~$K$, is $\Sp(K\otimes_KK)\cong\Sp(K)$ and therefore consists of a single point. This proves that $x'=x''$, as desired.

In particular, taking $Y'=X$, the projection $f_X:X\times_Y X\to X$ is injective, and hence the diagonal morphism $\Delta:X\to X\times_Y X$ is surjective (since $f_X\circ \Delta$ is the identity on $X$). On the other hand $\Delta$ is an open immersion, because $f$ is \'etale. It follows that $\Delta$ and hence $f_X$ are isomorphisms. On the other hand $f$ is flat by \'etaleness and even faithfully flat by surjectivity. Since being an isomorphism is local for the \'etale topology, and $f_X$ is an isomorphism, it follows that $f$ is an isomorphism, as desired.
\end{Proof}

\medskip
Now we look at the situation near the standard boundary component. 

\begin{Def}\label{def:NbhdOfInfty}
For any $n\in\BZ^{>0}$ and $R_n>0$ consider the $\Gamma_U$-invariant subset
$$I(n,R_n)\ :=\ \bigl\{\omega={\textstyle\binom{\omega_1}{\omega'}}\in\Omega^r 
\bigm| \omega'\in\Omega^{r-1}_n,\ d(\omega_1,\Finf^{r-1}\omega')\ge R_n\bigr\}.$$ 
An arbitrary $\Gamma_U$-invariant admissible open subset $\CN\subset \Omega^r$, such that for each $n>0$ there exists an $R_n>0$ with $I(n,R_n)\subset\CN$, will be called a \emph{neighbourhood of infinity}.
\end{Def}

Note that every subset of the form $I(n,R_n)$ is contained in $\Omega^r$ by construction; hence neighbourhoods of infinity exist and $\Omega^r$ is itself one.

\begin{Def}\label{def:quasi-uniform}
Any subset of\/ $\Cinf\times\Omega^{r-1}$ of the form
$$\CT\ =\ \bigcup_{n\ge1}B(0,r_n)\times\Omega^{r-1}_n$$
for numbers $r_n\in|\Cinf^\times|$ will be called a \emph{tubular neighbourhood  of} $\{0\}\times \Omega^{r-1}$, or just a {\em tubular neighbourhood} for the sake of brevity. The intersection of a tubular neighbourhood with $\Cinf^\times\times\Omega^{r-1}$ will be called a \emph{pierced tubular neighbourhood}.
\end{Def}

Any tubular neighbourhood is an admissible subset, because it is the union of affinoid sets of the form $B(0,\rho)\times \Omega^{r-1}_n$ for $\rho\in|\Cinf^\times|$ and the intersection of any two such affinoid sets is again of this form.
The same holds for pierced tubular neighbourhoods, but in this case we must use affinoids of the form $D(0,\sigma,\rho)\times \Omega^{r-1}_n$.

\medskip
Next recall that $e_{\Lambda'\omega'}(\omega_1)\not=0$ whenever $\omega_1\not\in\Lambda'\omega'$. In particular this holds for any $\omega = \binom{\omega_1}{\omega'}\in\Omega^r$, and so
\UseTheoremCounterForNextEquation
\begin{equation}\label{DefOfuParameter}
u_{\omega'}(\omega_1)\ :=\ \frac{1}{e_{\Lambda'\omega'}(\omega_1)}\ \in\ \Cinf^\times
\end{equation}
is well-defined for all $\omega = \binom{\omega_1}{\omega'}\in\Omega^r$.

\begin{Ex}\rm
Suppose that $A=\BF_q[t]$, $r=2$, $\Lambda=A^2$ and $\norm=\bar\pi$ is a period of the Carlitz module. Then for $\omega=\binom{\omega_1}{\norm}\in\Omega^2$ we have
\[u_{\omega'}(\omega_1) = \frac{1}{e_{\norm A}(\omega_1)} = \frac{1}{\bar\pi e_{A}(z)},\]
where $z=\omega_1/\norm\in\Cinf\smallsetminus\Finf$ is the usual parameter at infinity in the rank 2 literature (see, e.g., \cite{GekelerCoeff}). 
\end{Ex}


\begin{Thm}\label{Thm:Parameter}
\begin{enumerate}
\item[(a)] The morphism 
$$\textstyle\theta:\ 
\Gamma_U \backslash \Omega^r \longto\Cinf^\times\times\Omega^{r-1},\ \ 
\left[\binom{\omega_1}{\omega'}\right] \longmapsto \binom{u_{\omega'}(\omega_1)}{\omega'}$$
induces an isomorphism of rigid analytic spaces from $\Gamma_U \backslash \Omega^r$ to an admissible open subset of $\Cinf^\times\times\Omega^{r-1}$.
\item[(b)] For any neighbourhood of infinity $\CN\subset\Omega^r$, the image $\theta(\Gamma_U\backslash\CN)$ contains a pierced tubular neighbourhood.
\item[(c)] For any pierced tubular neighbourhood $\CT'\subset\Cinf^\times\times\Omega^{r-1}$ there is a neighbourhood of infinity $\CN\subset\Omega^r$ such that 
$\theta(\Gamma_U\backslash\CN)=\CT'$, and $\theta$ induces an isomorphism 
$$\Gamma_U\backslash\CN \stackrel{\sim}{\longto} \CT'.$$
\end{enumerate}
\end{Thm}

\begin{Proof}
Part (a) is a direct consequence of Proposition \ref{Prop:OmegaByGammaUIso}.
To prove (b) we must show that for any $n>0$ and $R_n>0$ there exists $r_n>0$ such that
\[B(0,r_n)'\times\Omega_n^{r-1}\ \subset\ \theta\bigl(\Gamma_U\backslash I(n,R_n) \bigr).\]
For this let $c_n$ be the constant from Proposition \ref{NewEarth} (c) and set $r_n := c_n^{-1}$. Consider any point $\binom{z}{\omega'}\in B(0,r_n)'\times\Omega_n^{r-1}$. As the map $e_{\Lambda'\omega'} : \Cinf\smallsetminus \Lambda'\omega' \to \Cinf^\times$ is surjective by Proposition \ref{exp1}, and $u_{\omega'}=e_{\Lambda'\omega'}^{-1}$ by definition, there exists a point $\omega_1\in\Cinf\setminus\Lambda'\omega'$ with $u_{\omega'}(\omega_1)=z$. Since $z\in B(0,r_n)'$, we then have $|e_{\Lambda'\omega'}(\omega_1)|\ge c_n$. By Proposition \ref{NewEarth} (c) we thus have $d(\omega_1,\Finf^{r-1}\omega')\ge R_n$, and so $\binom{\omega_1}{\omega'} \in I(n,R_n)$. Therefore $\binom{z}{\omega'} = \theta([\binom{\omega_1}{\omega'}]) \in \theta(\Gamma_U\backslash I(n,R_n))$, proving (b).

To prove (c) we must show that for any $n>0$ and $r_n>0$ there exists $R_n>0$ such that
$$\theta\bigl(\Gamma_U\backslash I(n,R_n)\bigr)
\ \subset\ B(0,r_n)'\times\Omega_n^{r-1}.$$
For this set $R_n:=r_n^{-1}$ and consider any point $\binom{\omega_1}{\omega'}\in I(n,R_n)$. Then by Proposition \ref{NewEarth} (d) we have
$|e_{\Lambda'\omega'}(\omega_1)|\ge d(\omega_1,F_\infty^{r-1}\omega') \ge R_n$ and hence $|u_{\omega'}(\omega_1)|\le r_n$. Therefore $\theta([\binom{\omega_1}{\omega'}]) \in B(0,r_n)'\times\Omega_n^{r-1}$, as desired. The isomorphy $\Gamma_U\backslash\CN \stackrel{\sim}{\longto} \CT'$ then follows from (a).
\end{Proof}

\section{Expansion at infinity}
\label{Sec:HolomorphyB}

In this section we show that every $\Gamma_U$-invariant holomorphic function 
admits a Laurent series expansion in $u_{\omega'}(\omega_1)$ which converges near infinity. 
As usual, we measure the size of a holomorphic function 
$g:\Omega^{r-1}_n\to\Cinf$ by the supremum norm
$$\|g\|_n:=\sup\{|g(\omega')|:\omega'\in\Omega^{r-1}_n\}.$$
Note that any rational function is bounded outside of a neighbourhood of its poles. In particular, a rational function with no poles on $\Omega^r$ is bounded on $\Omega^r_n$. Since $g$ is a \emph{uniform} limit of rational functions on $\Omega^r_n$, the supremum defined above will always be finite.

\begin{Lem}\label{Lem:BoundsOnCoefficients}
Let $n\in\BZ^{>0}$ and $\rho\in|\Cinf|$. Any holomorphic function $f:B(0,\rho)'\times\Omega^{r-1}_n\to\Cinf$ has a unique Laurent series expansion
\UseTheoremCounterForNextEquation
\begin{equation}\label{eq:LaurentExpansion}
f(z,\omega')=\sum_{k\in\BZ}f_k(\omega') z^k,
\end{equation}
which converges uniformly on every affinoid subset, where the functions $f_k:\Omega^{r-1}_n\to\Cinf$ are holomorphic and satisfy the conditions
$$ \limsup_{k\to\infty} \sqrt[k]{\|f_k\|_n} \ \le\ \rho^{-1}
\quad\text{and}\quad 
\lim_{k\to-\infty} \sqrt[-k]{\|f_k\|_n}\ =\ 0.$$
\end{Lem}

\begin{Proof}
Write $\rho = q^a$ with $a\in\BQ$. Then the punctured disk $B(0,\rho)'$ is the union of the affinoid annuli 
$$D(0,\sigma,\rho)\ =\ \{z\in\Cinf\mid \sigma\le|z|\le \rho\}
\ =\ \Spm\Cinf\langle \tfrac{X}{\pi^a},\tfrac{\pi^b}{X}\rangle$$
for all $\sigma = q^b<\rho$ with $b\in\BQ$. 
Since $\Omega_n^{r-1}$ is also affinoid, say $\Omega^{r-1}_n=\Spm A_n^{r-1}$, the product is affinoid and more precisely
$$D(0,\sigma,\rho)\times \Omega^{r-1}_n
\ =\ \Spm A_n^{r-1}\langle \tfrac{X}{\pi^a},\tfrac{\pi^b}{X}\rangle.$$
Thus the restriction of $f$ to $D(0,\sigma,\rho)\times \Omega^{r-1}_n$ has a uniformly convergent expansion of the form \eqref{eq:LaurentExpansion} with unique holomorphic functions $f_k:\Omega^{r-1}_n\to\Cinf$ that satisfy
$$\limsup_{k\to\infty} \sqrt[k]{\|f_k\|_n} \ \le\ \rho^{-1}
\quad\text{and}\quad 
\limsup_{k\to-\infty} \sqrt[-k]{\|f_k\|_n}\ \le\ \sigma.$$
By uniqueness, the functions $f_k$ are independent of~$\sigma$, so the proposition follows by letting $\sigma$ go to~$0$.
\end{Proof}

\begin{Lem}\label{Prop:LaurentExp}
For any pierced tubular neighbourhood $\CT'\subset \Cinf^\times\times\Omega^{r-1}$, any holomorphic function $f:\CT'\to\Cinf$ has a unique Laurent series expansion
$$f(z,\omega')=\sum_{k\in\BZ}f_k(\omega') z^k$$
with holomorphic functions $f_k:\Omega^{r-1}\to\Cinf$, which converges uniformly on every affinoid subset of~$\CT'$.
\end{Lem}

\begin{Proof}
Suppose that $\CT'=\bigcup_{n\ge1}B(0,r_n)'\times\Omega^{r-1}_n$ with $r_n\in|\Cinf^\times|$.
By Lemma~\ref{Lem:BoundsOnCoefficients}, for any $n\ge1$ the restriction of $f$ to $B(0,r_n)'\times\Omega^{r-1}_n$ admits a unique Laurent series expansion $\sum_{k\in\BZ} f^{(n)}_k z^k$ with holomorphic functions $f^{(n)}_k:\Omega^{r-1}_n\to\Cinf$ which converges uniformly on every affinoid subset. 
For any $n>m\ge1$, the uniqueness in Lemma~\ref{Lem:BoundsOnCoefficients} for the restriction of $f$ to $B(0,\min\{r_m,r_n\})'\times\Omega^{r-1}_m$ implies that $f^{(n)}_k|\Omega^{r-1}_m = f^{(m)}_k$. By the sheaf property for admissible coverings, there are therefore unique holomorphic functions $f_k:\Omega^{r-1}\to\Cinf$ with $f_k|\Omega^{r-1}_n = f^{(n)}_k$ for all~$n$, and they satisfy the desired conditions.
\end{Proof}

\begin{Prop}\label{ModFormsLaurentExpansion}
For any $\Gamma_U$-invariant holomorphic function $f:\Omega^r \to \Cinf$ there exist unique holomorphic functions $f_n: \Omega^{r-1} \to \Cinf$, such that the series
$$\sum_{n\in\BZ}f_n(\omega')\cdot u_{\omega'}(\omega_1)^n$$
converges to $f(\binom{\omega_1}{\omega'})$ on some neighbourhood of infinity, and uniformly on every affinoid subset thereof.
\end{Prop}

\begin{Proof}
Being $\Gamma_U$-invariant $f$ corresponds to a function $\bar f:\OmegaModGammaU\to\Cinf$. By Theorem \ref{Thm:Parameter} (c) the function $\bar f\circ\theta^{-1}$ then induces a holomorphic function on a pierced tubular neighbourhood $\CT'\subset\Cinf^\times \times\Omega^{r-1}$, where $\CT'=\theta(\Gamma_U\backslash\CN)\subset\Cinf^\times \times\Omega^{r-1}$
for a neighbourhood of infinity $\CN\subset\Omega^r$. By Lemma \ref{Prop:LaurentExp} that function has a unique expansion of the form
$${\textstyle\bar f\circ\theta^{-1}(\binom{z}{\omega'})}
\ =\ \sum_{n\in\BZ}f_n(\omega')z^n.$$
By the definition of $\theta$ this yields a unique expansion 
\[{\textstyle f(\binom{\omega_1}{\omega'})}\ =\ 
\sum_{n\in\BZ} f_n(\omega')\cdot u_{\omega'}(\omega_1)^n\]
on $\CN$, which again converges uniformly on every affinoid subset, as desired.
\end{Proof}

\begin{Rem}\label{ExpansionDoesNotConvergeEverywhere}\rm
The series in Proposition \ref{ModFormsLaurentExpansion} does not necessarily converge on all of~$\Omega^r$. For example, in \cite[Corollary 2.2]{GekelerGrowth}, Gekeler shows that the $u$-expansion of the rank 2 Drinfeld discriminant function has the radius of convergence $q^{-1/(q-1)}$ only. This is in contrast with the classical case, where the $q$-expansion of a modular form converges on the entire upper half plane.  
\end{Rem}

Any weak modular form for the group $\Gamma$ is a $\Gamma_U$-invariant function; hence it possesses a $u$-expansion as in Proposition \ref{ModFormsLaurentExpansion}.
Our next aim is to study its behaviour under conjugation by the ``stabiliser of the standard boundary component''. For this consider the algebraic subgroups 
\UseTheoremCounterForNextEquation
\begin{equation}\label{DefPL}
P\ :=\ \left(\begin{array}{c|c}
* & *\ \ldots\ * \\ \hline
\begin{array}{c} 0 \\[-5pt] \vdots \\[-4pt] 0 \end{array}
& \begin{array}{c} *\ \ldots\ * \\[-5pt] \vdots \ \qquad \vdots \\[-4pt]  *\ \ldots\ * \end{array}
\end{array}\right)_,
\quad \text{and}\quad
M\ :=\ 
\left(\begin{array}{c|c}
* & 0\ \ldots\ 0 \\ \hline
\begin{array}{c} 0 \\[-5pt] \vdots \\[-4pt] 0 \end{array}
& \begin{array}{c} *\ \ldots\ * \\[-5pt] \vdots \ \qquad \vdots \\[-4pt]  *\ \ldots\ * \end{array}
\end{array}\right)_.
\end{equation}
of $\GL_{r,F}$, so that $P=U\rtimes M$ is parabolic with unipotent radical $U$ and Levi subgroup~$M$. 

\begin{Lem}\label{ActionOfL}
Consider any element of the form $\gamma=\binom{\alpha\ \kern2pt 0\kern1pt}{0\ \gamma'} \in M(F)$ with $\alpha\in F^\times$ and $\gamma\in\GL_{r-1}(F)$ and any point $\omega = \binom{\omega_1}{\omega'}\in\Omega^r$. Then:
\begin{enumerate}
\item[(a)] $\eta := j(\gamma,\omega)=j(\gamma',\omega')$ and $\gamma(\omega) = \binom{\eta^{-1}\alpha\omega_1}{\gamma'(\omega')}$. 
\item[(b)] $\Lambda'_\gamma := \iota^{-1}(\gamma^{-1}\Gamma_U\gamma) =\alpha^{-1}\Lambda'\gamma'$.
\item[(c)] $u_{\gamma,\omega'}(\omega_1) := e_{\Lambda'_\gamma\omega'}(\omega_1)^{-1} = \eta^{-1}\alpha\cdot u_{\gamma'(\omega')}(\eta^{-1}\alpha\omega_1)$.
\item[(d)] There exist constants $k\ge0$ and $c_4>0$ such that for all $n>0$ and $R>0$ we have
$$\gamma\left(I(n,R)\right)\ \subset\ I(n+k,c_4R).$$
\item[(e)] For any neighbourhood of infinity $\CN\subset\Omega^r$ the subset $\gamma^{-1}(\CN)$ is also a neighbourhood of infinity.
\item[(f)] For any $\Gamma_U$-invariant holomorphic function $f:\Omega^r \to \Cinf$ with the expansion in Proposition \ref{ModFormsLaurentExpansion} on~$\CN$ and any integers $k$ and $m$ we have the following expansion on~$\gamma^{-1}(\CN)$:
$${\textstyle (f|_{k,m}\gamma)(\binom{\omega_1}{\omega'})}\ =\ 
\sum_{n\in\BZ}\alpha^{m-n}(f_n|_{k-n,m}\gamma')(\omega')\cdot u_{\gamma,\omega'}(\omega_1)^n.$$
\end{enumerate}
\end{Lem}

\begin{Proof}
Assertion (a) follows directly from the definitions (\ref{DefOfj}) and (\ref{DefOfActionOnOmega}), with $\gamma'(\omega') = \eta^{-1}\gamma'\omega'$.
Assertion (b) follows by direct calculation from the definition (\ref{eq:iota}) of~$\iota$. Using (b) and Proposition \ref{exp2} (b) we deduce that
$$e_{\Lambda'_\gamma\omega'}(\omega_1)
\ =\ e_{\alpha^{-1}\Lambda'\gamma'\omega'}(\omega_1)
\ =\ e_{\alpha^{-1}\eta\Lambda'\cdot\gamma'(\omega')}(\omega_1)
\ =\ \alpha^{-1}\eta\cdot e_{\Lambda'\cdot\gamma'(\omega')}(\eta^{-1}\alpha\omega_1)$$
Taking reciprocals thus shows (c).

To prove (d) consider any $n>0$ and $\omega'\in\Omega^{r-1}_n$. Then by definition (\ref{Omega_n}) and Lemma \ref{Lem:BoundingAbsoluteValues} (c), both with $r-1$ in place of~$r$, we have $h(\omega')\ge |\pi|^n$ and $h(\gamma'(\omega')) \ge c_3 h(\omega')$ for some constant $c_3$ depending only on $\gamma'$. Together we deduce that $h(\gamma'(\omega')) \ge |\pi|^{n+k}$ for some $k\ge0$ depending only on~$\gamma'$. By the definition (\ref{Omega_n}) again this means that $\gamma'(\omega')\in\Omega^{r-1}_{n+k}$.
Next, by Lemmas \ref{Lem:BoundingAbsoluteValues} (a) and \ref{Lem:Omeganbounds}, again with $r-1$ in place of~$r$, we have $|\eta| = |j(\gamma',\omega')|\le |\omega'|c_1^{-1} \leq |\xi||\pi|^{-n}c_1^{-1}$ for another constant $c_1$ depending only on~$\gamma'$. 
Note also that, since $\gamma'(\omega') = \eta^{-1}\gamma'\omega'$, the associated $F_\infty$-vector space is $\Finf^{r-1}\gamma'(\omega') = \eta^{-1}\Finf^{r-1}\omega'$. For any $\omega_1\in\Cinf$ we therefore have
\begin{eqnarray*}
d(\eta^{-1}\alpha\omega_1,\Finf^{r-1}\gamma'(\omega'))
\!\!&=&\!\! d(\eta^{-1}\alpha\omega_1,\eta^{-1}\alpha\Finf^{r-1}\omega') \\
\!\!&=&\!\! |\eta^{-1}\alpha|\cdot d(\omega_1,\Finf^{r-1}\omega')
\ \ge\ |\alpha\pi^n\xi^{-1}|c_1\cdot d(\omega_1,\Finf^{r-1}\omega').
\end{eqnarray*}
In view of Definition \ref{def:NbhdOfInfty} this implies (d) with $c_4 := |\alpha\pi^n\xi^{-1}|c_1$.

To deduce (e) choose $R_n>0$ such that $\bigcup_{n>0}I(n,R_n)\subset\CN$. Then (d) implies that 
$$\textstyle \gamma\bigl(\bigcup_{n>0}I(n,c_4^{-1}R_{n+k})\bigr)
\ \subset\ \bigcup_{n>0}I(n+k,R_{n+k})
\ \subset\ \CN,$$
and hence $\bigcup_{n>0}I(n,c_4^{-1}R_{n+k}) \subset \gamma^{-1}(\CN)$. Thus $\gamma^{-1}(\CN)$ is a neighbourhood of infinity, proving (e).

Finally, using the definition (\ref{DefOfActionOnF}), for any ${\textstyle\binom{\omega_1}{\omega'}}\in\gamma^{-1}(\CN)$ we can now calculate
\begin{eqnarray*}
{\textstyle (f|_{k,m}\gamma)(\binom{\omega_1}{\omega'})}
&\stackrel{\rm (a)}{=}& \eta^{-k}(\det \gamma)^{m}f({\textstyle\binom{\eta^{-1}\alpha\omega_1}{\gamma'(\omega')}})\\
&\stackrel{\ref{ModFormsLaurentExpansion}}{=}& \eta^{-k}(\det\gamma)^{m} \cdot\sum_{n\in\BZ}f_n(\gamma'(\omega'))\cdot u_{\gamma'(\omega')}(\eta^{-1}\alpha\omega_1)^n\\
&\stackrel{\rm (c)}{=}& \eta^{-k}(\alpha\det\gamma')^{m}\cdot\sum_{n\in\BZ}f_n(\gamma'(\omega'))\cdot \bigl(\alpha^{-1}\eta u_{\gamma,\omega'}(\omega_1))^n\\
&=& \sum_{n\in\BZ}\alpha^{m-n}\cdot\eta^{n-k}(\det\gamma')^{m}f_n(\gamma'(\omega'))\cdot u_{\gamma,\omega'}(\omega_1)^n\\
&=&\kern-5pt \sum_{n\in\BZ}\alpha^{m-n}\cdot(f_n|_{k-n,m}\gamma')(\omega')\cdot u_{\gamma,\omega'}(\omega_1)^n,
\end{eqnarray*}
proving (f).
\end{Proof}

\medskip
For a first application 
consider the subgroup
\UseTheoremCounterForNextEquation
\begin{equation}\label{GammaLeviDef}
\textstyle \Gamma_M\ :=\ 
\bigl\{ \gamma'\in\GL_{r-1}(F)
\bigm| \binom{1\ \kern2pt 0\kern1pt}{0\ \gamma'} \in \Gamma\cap M(F)
\bigr\}.
\end{equation}

\begin{Thm}\label{Thm:CoefficientsAreModularForms}
Let $f$ be a weak modular form of weight $k$ and type $m$ for the group $\Gamma$, and let $f_n$ be its coefficients in the $u$-expansion from Proposition \ref{ModFormsLaurentExpansion}.
Then, for each $n\in\BZ$, the function $f_n$ is a weak modular form of weight $k-n$ and type $m$ for the group $\Gamma_M<\GL_{r-1}(F)$.
\end{Thm}

\begin{Proof}
Consider any $\gamma'\in\Gamma_M$ and set $\gamma:=\binom{1\ \kern2pt 0\kern1pt}{0\ \gamma'}$, so that $\alpha=1$ in the notation of Lemma \ref{ActionOfL}. Since the subgroup $\Gamma_U$ is normalised by~$\gamma$, Lemma \ref{ActionOfL} (b) implies that $\Lambda'_{\gamma}=\Lambda'$ and hence $u_{\gamma,\omega'}(\omega_1) = u_{\omega'}(\omega_1)$.
Let $\CN$ be a neighbourhood of infinity on which the expansion from Proposition \ref{ModFormsLaurentExpansion} converges. Then by Lemma \ref{ActionOfL} (e) the intersection $\CN\cap\gamma^{-1}(\CN)$ is another neighbourhood of infinity.
For any $\omega = \binom{\omega_1}{\omega'}\in \CN\cap\gamma^{-1}(\CN)$ we can compare the expansions of $f(\omega)=(f|_{k,m}\gamma)(\omega)$ from 
Proposition \ref{ModFormsLaurentExpansion} and \ref{ActionOfL} (f). Since $\alpha=1$, by the uniqueness part of Proposition \ref{ModFormsLaurentExpansion} we conclude that $f_n = f_n|_{k-n,m}\gamma'$ for all $n\in\BZ$, proving the theorem.
\end{Proof}

\begin{Cor}\label{Cor:SomeCoeffsVanish}
Let $f$ be a weak modular form of weight $k$ and type $m$ for the group $\Gamma$. Then for any $n\in\BZ$, the coefficient $f_n$ in the $u$-expansion from Proposition \ref{ModFormsLaurentExpansion} is identically zero unless
$$n \equiv k-(r-1)m \text{\ \ modulo\ \ } |\Gamma_M\cap\{\text{scalars}\}|.$$
\end{Cor}

\begin{Proof}
Combine Theorem \ref{Thm:CoefficientsAreModularForms} with (\ref{TypeZero}) for $r-1$ in place of~$r$.
\end{Proof}

\begin{Lem}\label{ActionOfU}
Consider any element of the form $\gamma=\binom{1\kern7pt \beta\kern10pt}{0\ \Id_{r-1}} \in U(F)$ for some row vector $\beta\in F^{r-1}$ and any point $\omega = \binom{\omega_1}{\omega'}\in\Omega^r$. Then:
\begin{enumerate}
\item[(a)] $j(\gamma,\omega)=\det(\gamma)=1$ and $\gamma(\omega) = \binom{\omega_1+\beta\omega'}{\omega'}$. 
\item[(b)] For any neighbourhood of infinity $\CN\subset\Omega^r$ the subset 
$$\CN'\ :=\ \bigl\{
{\textstyle\binom{\omega_1}{\omega'}} \in \gamma^{-1}(\CN) \bigm|
|e_{\Lambda'\omega'}(\beta\omega') \cdot u_{\omega'}(\omega_1)|<1\bigr\}$$
is also a neighbourhood of infinity.
\item[(c)] For any $\Gamma_U$-invariant holomorphic function $f:\Omega^r \to \Cinf$ with the expansion in Proposition \ref{ModFormsLaurentExpansion} on~$\CN$ and any integers $k$ and $m$ we have the following expansion on~$\CN'$:
$${\textstyle (f|_{k,m}\gamma)(\binom{\omega_1}{\omega'})}
\ =\ \sum_{n\in\BZ} \;\Bigl(\;\sum_{k\ge0} {\textstyle\binom{k-n}{k}}\cdot f_{n-k}(\omega')\cdot e_{\Lambda'\omega'}(\beta\omega')^{k} \Bigr) \cdot u_{\omega'}(\omega_1)^n.$$
\end{enumerate}
\end{Lem}

\begin{Proof}
Assertion (a) follows directly from the definitions (\ref{DefOfj}) and (\ref{DefOfActionOnOmega}). 

To prove (b) choose $R_n>0$ such that $\bigcup_{n>0}I(n,R_n)\subset\CN$. Since $\beta\omega'\in \Finf^{r-1}\omega'$, we have $d(\omega_1+\beta\omega',\Finf^{r-1}\omega') = d(\omega_1,\Finf^{r-1}\omega')$ and therefore $\gamma^{-1}(I(n,R_n))=I(n,R_n)$ by Definition \ref{def:NbhdOfInfty}. 
On the other hand we have $d(\beta\omega',F_\infty^{r-1}\omega')=0$; applying Proposition \ref{NewEarth} (c) thus yields constants $c_n>0$, such that $|e_{\Lambda'\omega'}(\beta\omega')|<c_n$ for any $\omega'\in\Omega^{r-1}_n$. 
By Proposition \ref{NewEarth} (d) and Definition \ref{def:NbhdOfInfty}, for any $\binom{\omega_1}{\omega'}\in I(n,c_n)$ we therefore have 
$$|e_{\Lambda'\omega'}(\omega_1)|\ \ge\ d(\omega_1,\Finf^{r-1}\omega')\ \ge\ c_n\ >\ |e_{\Lambda'\omega'}(\beta\omega')|.$$
By the definition of $u_{\omega'}(\omega_1)$ this implies that $|e_{\Lambda'\omega'}(\beta\omega') \cdot u_{\omega'}(\omega_1)|<1$.
Together this shows that $I(n,\max\{R_n,c_n\})\subset\CN'$.
Varying $n$ we conclude that $\CN'$ is a neighbourhood of infinity, proving (b).

Next, by (a) and the definition (\ref{DefOfActionOnF}), the expansion from Proposition \ref{ModFormsLaurentExpansion} yields
$${\textstyle (f|_{k,m}\gamma)(\binom{\omega_1}{\omega'})}
\ =\ f({\textstyle\binom{\omega_1+\beta\omega'}{\omega'}})
\ =\ \sum_{n\in\BZ}f_n(\omega')\cdot u_{\omega'}(\omega_1+\beta\omega')^n$$
Using the additivity of the exponential function we can rewrite
\begin{eqnarray*}
u_{\omega'}(\omega_1+\beta\omega')^n
&=& e_{\Lambda'\omega'}(\omega_1+\beta\omega')^{-n} \\
&=& \bigl(e_{\Lambda'\omega'}(\omega_1)+e_{\Lambda'\omega'}(\beta\omega')\bigr)^{-n} \\
&=& \bigl(1+e_{\Lambda'\omega'}(\beta\omega') u_{\omega'}(\omega_1)\bigr)^{-n} \cdot u_{\omega'}(\omega_1)^n.
\end{eqnarray*}
For $\binom{\omega_1}{\omega'} \in \CN'$ we have $|e_{\Lambda'\omega'}(\beta\omega') \cdot u_{\omega'}(\omega_1)|<1$, so we can plug the binomial series into the above expansion and rearrange terms, yielding
\begin{eqnarray*}
{\textstyle (f|_{k,m}\gamma)(\binom{\omega_1}{\omega'})}
&=& \sum_{n\in\BZ} f_n(\omega')\cdot \Bigl( \sum_{k\ge0} {\textstyle\binom{-n}{k}}\cdot e_{\Lambda'\omega'}(\beta\omega')^{k} \cdot u_{\omega'}(\omega_1)^{n+k} \Bigr) \\
&=& \sum_{n\in\BZ} \sum_{k\ge0} {\textstyle\binom{-n}{k}}\cdot f_n(\omega')\cdot e_{\Lambda'\omega'}(\beta\omega')^{k} \cdot u_{\omega'}(\omega_1)^{n+k} \\
&=& \sum_{n'\in\BZ} \;\Bigl(\;\sum_{k\ge0} {\textstyle\binom{k-n'}{k}}\cdot f_{n'-k}(\omega')\cdot e_{\Lambda'\omega'}(\beta\omega')^{k} \Bigr) \cdot u_{\omega'}(\omega_1)^{n'}
\end{eqnarray*}
with the substitution $n+k=n'$. Thus the stated expansion holds on~$\CN'$, proving (c).
\end{Proof}

\begin{Def}\label{Def:OrderAtInfinity}
Let $f : \Omega^r \to \Cinf$ be a $\Gamma_U$-invariant holomorphic function and let $f_n$ be its coefficients in the $u$-expansion from Proposition \ref{ModFormsLaurentExpansion}.
Then the \emph{order at infinity} of $f$  is
$$\ord_{\Gamma_U}(f)\ :=\ \inf\bigl\{n\in\BZ \bigm| f_n(\omega')\ne 0\text{ for some }\omega'\in\Omega^{r-1}\bigr\}
\ \in\ \BZ\cup\{\pm\infty\}.$$
The function $f$ is called \emph{meromorphic at infinity} if $\ord_{\Gamma_U}(f)>-\infty$, that is, if $f_n$ is identically zero for all $n\ll0$.
It is called \emph{holomorphic at infinity} if $\ord_{\Gamma_U}(f)\ge0$, that is, if $f_n$ is identically zero for all $n<0$.
\end{Def}

\begin{Prop}\label{Prop:OrderInvariantA}
Consider a $\Gamma_U$-invariant holomorphic function $f \colon \Omega^r \to \Cinf$ and an element $\gamma\in P(F)$. Then $f|_{k,m}\gamma$ is invariant under 
$\Gamma_{\gamma,U} := (\gamma^{-1}\Gamma\gamma)\cap U(F)$, and we have
$$\ord_{\Gamma_U}(f) = \ord_{\Gamma_{\gamma,U}}(f|_{k,m}\gamma).$$
In particular $f$ is meromorphic, respectively holomorphic at infinity if and only if $f|_{k,m}\gamma$ has the corresponding property.
\end{Prop}

\begin{Proof}
Since $P=U\rtimes M$, it suffices to prove this separately for elements of $M(F)$ and $U(F)$. In both cases the $\Gamma_{\gamma,U}$-invariance follows by direct calculation from the formula (\ref{fCocycle}). 
The rest follows from the expansion in Lemma \ref{ActionOfL} for $\gamma\in M(F)$, respectively by close inspection of the expansion in Lemma \ref{ActionOfU} for $\gamma\in U(F)$.
\end{Proof}

\begin{Prop}\label{Prop:Order}
Let $\Gamma_1<\Gamma$ and hence $\Gamma_{1,U}:=\Gamma_1\cap U(F) < \Gamma_U$ be subgroups of 
finite index. Then for any $\Gamma_U$-invariant holomorphic function $f$ we have
$$\ord_{\Gamma_{1,U}}(f) = \ord_{\Gamma_U}(f)\cdot[\Gamma_U : \Gamma_{1,U}].$$
In particular $f$ is meromorphic, respectively holomorphic at infinity with respect to $\Gamma_U$ 
if and only if it is so with respect to $\Gamma_{1,U}$.
\end{Prop}

\begin{Proof}
The lattice associated to $\Gamma_{1,U}$ is $\Lambda_1' := \iota^{-1}(\Gamma_{1,U}) \subset \Lambda' = \iota^{-1}(\Gamma_U)$, so that $[\Lambda':\Lambda_1'] = [\Gamma_U:\Gamma_{1,U}] = p^d$ for an integer $d\ge0$. 
For any $\omega'\in\Omega^{r-1}$ we then also have $[\Lambda'\omega':\Lambda_1'\omega'] = p^d$. Let $B$ be a set of representatives for $\Lambda'\setminus\Lambda_1'$ modulo $\Lambda_1'$. By \Proposition \ref{exp2} (a) we then have 
$$e_{\Lambda'\omega'}(\omega_1)\ =\ e_{\Lambda_1'\omega'}(\omega_1) \cdot
\prod_{\beta\in B} \left( 1 - \frac{e_{\Lambda_1'\omega'}(\omega_1)} {e_{\Lambda_1'\omega'}(\beta\omega')} \right)_.$$
Taking reciprocals, we can therefore express the expansion parameter $u_{\omega'}(\omega_1) := e_{\Lambda'\omega'}(\omega_1)^{-1}$ with respect to $\Lambda'$ in terms of the expansion parameter $u_{1,\omega'}(\omega_1):= e_{\Lambda_1'\omega'}(\omega_1)^{-1}$ with respect to $\Lambda_1'$ by the formula
$$u_{\omega'}(\omega_1)\ =\ 
u_{1,\omega'}(\omega_1)^{p^d} \cdot \prod_{\beta\in B} 
\frac{e_{\Lambda_1'\omega'}(\beta\omega')} {e_{\Lambda_1'\omega'}(\beta\omega')u_{1,\omega'}(\omega_1)-1} _.$$
The expansion from Proposition \ref{ModFormsLaurentExpansion} thus yields
$$f({\textstyle\binom{\omega_1}{\omega'}})\ =\ 
\sum_{n\in\BZ}f_n(\omega')\cdot u_{\omega'}(\omega_1)^n\ =\ 
\sum_{n\in\BZ}f_n(\omega')\cdot 
u_{1,\omega'}(\omega_1)^{np^d} \cdot \prod_{\beta\in B} 
\left(\frac{e_{\Lambda_1'\omega'}(\beta\omega')} {e_{\Lambda_1'\omega'}(\beta\omega')u_{1,\omega'}(\omega_1)-1} \right)^{n}$$
for all points $\binom{\omega_1}{\omega'}$ in some neighbourhood of infinity.
By Lemma \ref{ActionOfU} (b) with $\Gamma_{1,U}$ in place of $\Gamma_U$, for each $\beta\in B$ we have $\bigl|e_{\Lambda_1'\omega'}(\beta\omega')u_{1,\omega'}(\omega_1)\bigr|<1$ on some neighbourhood of infinity.
On the intersection of these neighbourhoods, we can plug the binomial series into the above expansion and rearrange terms. We conclude that the expansion with respect to $u_{\omega'}(\omega_1)$ has the first non-zero term $f_n(\omega')\cdot u_{\omega'}(\omega_1)^n$ if and only if the expansion with respect to $u_{1,\omega'}(\omega_1)$ has the first non-zero term
$$f_n(\omega')\cdot u_{1,\omega'}(\omega_1)^{np^d} \cdot \prod_{\beta\in B} 
\left(-e_{\Lambda_1'\omega'}(\beta\omega') \right)^{n}.$$
Then $\ord_{\Gamma_{1,U}}(f) = np^d = \ord_{\Gamma_U}(f)\cdot[\Gamma_U : \Gamma_{1,U}]$, and the proposition follows.
\end{Proof}

\section{Modular forms}
\label{Sec:AMF}

Now we impose holomorphy conditions at all boundary components, not just the standard one. We achieve this by conjugating the standard boundary component by arbitrary elements $\delta\in\GL_r(F)$. 
Recall from \Proposition \ref{Prop:WConj} that for any weak modular form $f$ of weight $k$ and type $m$ for~$\Gamma$, and for any $\delta\in\GL_r(F)$, the function $f|_{k,m}\delta$ is a weak modular form of weight $k$ and type $m$ for the arithmetic subgroup $\delta^{-1}\Gamma\delta$. Determining the behaviour of $f$ at all boundary components is equivalent to determining the behaviour of all conjugates $f|_{k,m}\delta$ at the standard boundary component.


\begin{Def}\label{Def:ModularForms}
Let $f$ be a weak modular form of weight $k$ and type $m$ for $\Gamma$.
\begin{enumerate}
\item[(a)] If $\ord_{(\delta^{-1}\Gamma\delta)\cap U(F)}(f|_{k,m}\delta)\geq 0$ for all $\delta\in\GL_r(F)$, we call $f$ a {\em modular form}.
\item[(b)] If $\ord_{(\delta^{-1}\Gamma\delta)\cap U(F)}(f|_{k,m}\delta)\geq 1$ for all $\delta\in\GL_r(F)$, we call $f$ a {\em cusp form}.
\end{enumerate}
In particular, a \emph{modular form} is a weak modular form $f$ such that $f|_{k,m}\delta$ is holomorphic at infinity for all $\delta\in\GL_r(F)$.
The space of these functions is denoted by $\CM_{k,m}(\Gamma)$.
The space of cusp forms is denoted by $\CS_{k,m}(\Gamma)$.
As with weak modular forms, we abbreviate $\CM_k(\Gamma) := \CM_{k,0}(\Gamma)$ and $\CS_k(\Gamma) := \CS_{k,0}(\Gamma)$.
\end{Def}

It may seem extravagant to impose conditions for infinitely many $\delta$. However, the next two facts show that for fixed~$\Gamma$, we only need to check these conditions for $\delta$ in a fixed finite set.

\begin{Prop}\label{Prop:OrderInvariantB}
The numbers in Definition \ref{Def:ModularForms} depend only on the double coset $\Gamma\delta P(F)$.
\end{Prop}

\begin{Proof}
Since $f$ is a weak modular form of weight $k$ and type $m$ for $\Gamma$, for any $\delta'=\gamma'\delta\gamma$ with $\gamma'\in\Gamma$ and $\gamma\in P(F)$ we have $f|_{k,m}\delta' = (f|_{k,m}\delta)|_{k,m}\gamma$ and hence $\ord_{(\delta^{\prime-1}\Gamma\delta')\cap U(F)}(f|_{k,m}\delta')
 = \ord_{(\delta^{-1}\Gamma\delta)\cap U(F)}(f|_{k,m}\delta)$ by Proposition \ref{Prop:OrderInvariantA}.
\end{Proof}

\begin{Prop}\label{Prop:cosets}
The double coset space $\Gamma\backslash\GL_r(F)/P(F)$ is finite. More precisely, let $\Cl(A)$ denote the class group of $A$. Then:
\begin{enumerate}
\item[(a)] $\GL_r(A)\backslash\GL_r(F)/P(F)$ is in bijection with $\Cl(A)$.
\item[(b)] For any arithmetic subgroup $\Gamma<\GL_r(F)$, the set $\Gamma\backslash\GL_r(F)/P(F)$ has cardinality at most $|\Cl(A)|\cdot[\GL_r(A) : \GL_r(A)\cap\Gamma]$.
\item[(c)] If $\Gamma<\GL_r(A)$ then the double cosets of $\Gamma\backslash\GL_r(F)/P(F)$ can be represented by elements of $\GL_r(A)$ if and only if $\Cl(A)=\{1\}$. \label{Prop:cosetsforGLrA}
\end{enumerate}
\end{Prop}

\begin{Proof}
By the orbit-stabiliser theorem the set $\GL_r(F)/P(F)$ is in bijection with the set of one-dimensional subspaces of $F^r$ and hence with $\BP^{r-1}(F)$. This bijection is equivariant under the left action of $\GL_r(F)$. To prove (a) it thus suffices to find a bijection between $\GL_r(A)\backslash\BP^{r-1}(F)$ and $\Cl(A)$.

For this we associate to any column vector $x=(x_i)_i\in F^r\setminus\{0\}$ the fractional ideal $I(x) := \sum_i Ax_i \subset F$. This ideal depends only on the $\GL_r(A)$-orbit of~$x$, and its ideal class depends only on the corresponding point of $\BP^{r-1}(F)$. Together we therefore obtain a well-defined map $\GL_r(A)\backslash\BP^{r-1}(F) \to \Cl(A)$. This map is surjective, because $r\geq 2$ and every ideal of a Dedekind domain can be generated by $2$ elements. We claim that it is also injective. 

To see this we view $A^r$ as a space of row vectors, so that right multiplication by $x$ determines a surjective homomorphism of $A$-modules $p_x: A^r\to I(x)$. Since $I(x)$ is a projective $A$-module, the associated short exact sequence $0\to\ker(p_x)\to A^r\to I(x)\to 0$ splits. Moreover, since the isomorphism class of a finitely generated projective $A$-module depends only on its rank and its highest exterior power, the isomorphism class of $\ker(p_x)$ is determined by that of $I(x)$. 

Suppose now that two vectors $x$, $y\in F^r\setminus\{0\}$ correspond to the same ideal class. Then $I(y)=u\cdot I(x)$ for some $u\in F^\times$, and by the preceding remarks there exists an isomorphism of $A$-modules $f: \ker(p_x) \to \ker(p_y)$. Combining these via suitable splittings we find an isomorphism of $A$-modules $g: A^r\to A^r$ making the following diagram commute:
$$\xymatrix{
0\ar[r] & \ker(p_x) \ar[d]_{f}^\wr \ar[r] & A^r\ar[d]_{g}^\wr \ar[r]^{p_x} & I(x) \ar[d]_{u\cdot}^\wr \ar[r] & 0 \\
0\ar[r] & \ker(p_y) \ar[r] & A^r\ar[r]^{p_y} & I(y) \ar[r] & 0\rlap{.} \\}$$
Writing $g$ as right multiplication by a matrix $\gamma\in\GL_r(A)$, the commutativity on the right hand side then means that $a\gamma y=axu$ for all $a\in A^r$. Thus $\gamma y=xu$ for some $\gamma\in\GL_r(A)$ and $u\in F^\times$, which is precisely the desired injectivity.

This finishes the proof of (a). Parts (b) and (c) are direct consequences of (a).
\end{Proof}

\begin{Cor}\label{Cor:GLrACase}
Suppose that $\Gamma = \GL_r(A)$ for a principal ideal domain~$A$. Then:
\begin{enumerate}
\item[(a)] The condition in Definition \ref{Def:ModularForms} is independent of~$\delta$.
\item[(b)] If $m\not\equiv 0\ {\rm mod}\ (q-1)$, any modular form of weight $k$ and type $m$ for $\Gamma$ is a cusp form.
\end{enumerate}
\end{Cor}

\begin{Proof}
Part (a) follows from Propositions \ref{Prop:OrderInvariantB} and \ref{Prop:cosets} (a). To prove (b) let $f$ be a modular form of weight $k$ and type $m$ for $\Gamma$, and let $f_n$ be its coefficients in the $u$-expansion from Proposition \ref{ModFormsLaurentExpansion}, which are weak modular forms for the group $\Gamma_M=\GL_{r-1}(A)$. By assumption we then have $f_n=0$ for all $n<0$. If $f$ is not a cusp form, then $f_0$ is not identically zero, so 
Corollary \ref{Cor:SomeCoeffsVanish} implies that $k \equiv (r-1)m$ modulo $|\Gamma_M\cap\{\text{scalars}\}|=q-1$. But then $f$ itself is also not identically zero, so (\ref{TypeZero}) gives $k \equiv rm$ modulo $|\Gamma\cap\{\hbox{scalars}\}|=q-1$. Both congruences together imply that $m \equiv 0$ modulo $(q-1)$, contrary to the assumption.
\end{Proof}

\begin{Rem}\label{Thm:CoefficientsAreModularFormsRem}
\rm By Theorem \ref{Thm:CoefficientsAreModularForms} the coefficient $f_n$ of the $u$-expansion of a modular form $f$ is a weak modular forms of weight $k-n$ for a subgroup $\Gamma_M<\GL_{r-1}(F)$. In contrast to the case of modular forms in characteristic zero, the weight $k-n$ here goes to~$-\infty$ for $n\to\infty$. In Theorem \ref{FiniteDimension} (b) of Part II
we will see that any modular forms of weight $<0$ for $\Gamma_M$ must be zero if $r-1\ge2$. Thus for $r\ge3$, the coefficients $f_n$ cannot all be modular forms if they are non-zero.
\end{Rem}

\begin{Prop}\label{Prop:Conj}
For any $\delta\in\GL_r(F)$ we have $f\in \CM_{k,m}(\Gamma)$ if and only if 
$f|_{k,m}\delta \in \CM_{k,m}(\delta^{-1}\Gamma\delta)$.
\end{Prop}

\begin{Proof}
Direct consequence of 
\Proposition \ref{Prop:WConj} and the formula (\ref{fCocycle}).
\end{Proof}

\medskip
In particular, whenever $\Gamma_1\triangleleft\Gamma$ is a normal subgroup of finite index, the map $f\mapsto f|_{k,m}\gamma$ for all $\gamma\in\Gamma$ defines a right action of $\Gamma$ on $\CM_{k,m}(\Gamma_1)$. As a direct consequence of \Definition \ref{Def:ModularForms} and Proposition \ref{Prop:Order} the subspace of invariants is then
\UseTheoremCounterForNextEquation
\begin{equation}\label{ModInvs}
\CM_{k,m}(\Gamma_1)^{\Gamma} = \CM_{k,m}(\Gamma).
\end{equation}
Moreover, (\ref{TypeCongruence}) and (\ref{TypeZero}) imply that
\UseTheoremCounterForNextEquation
\begin{eqnarray}\label{TypeCongruenceM}
\CM_{k,m}(\Gamma) \!\!&=&\!\! \CM_{k,m'}(\Gamma)
\hbox{\ whenever $m\equiv m'$ modulo $|\det(\Gamma)|$, and}\\[3pt]
\UseTheoremCounterForNextEquation
\CM_{k,m}(\Gamma) \!\!&=&\!\! 0
\hbox{\ unless $k \equiv rm$ modulo $\bigl|\Gamma\cap\{\hbox{scalars}\}\bigr|$.}
\end{eqnarray}
As a direct consequence of the definitions we also have
\UseTheoremCounterForNextEquation
\begin{equation}\label{ModProduct}
\CM_{k,m}(\Gamma)\cdot\CM_{k',m'}(\Gamma) \subset \CM_{k+k',m+m'}(\Gamma)
\end{equation}
for all $k,k',m,m'$. In particular we can form the \emph{graded ring of modular forms}
\UseTheoremCounterForNextEquation
\begin{equation}\label{ModRing}
\CM_{*}(\Gamma) := \bigoplus_{k\ge0}\CM_{k}(\Gamma).
\end{equation}



\begin{center}
\rule{8cm}{0.01cm}
\end{center}

\begin{minipage}[t]{5cm}{\small
Department of Mathematical Sciences \\
University of Stellenbosch \\
Stellenbosch, 7600 \\
South Africa \\
djbasson@sun.ac.za
}
\end{minipage}\hfill
\begin{minipage}[t]{5cm}{\small
School of Mathematical and Physical Sciences \\
University of Newcastle \\
Callaghan, 2308 \\
Australia \\
florian.breuer@newcastle.edu.au\\
 {\em and} \\
Department of Mathematical Sciences \\
University of Stellenbosch \\
Stellenbosch, 7600 \\
South Africa \\}
\end{minipage}\hfill
%
%
\begin{minipage}[t]{5cm}{\small
Department of Mathematics \\
ETH Z\"urich\\
8092 Z\"urich\\
Switzerland \\
pink@math.ethz.ch}
\end{minipage}

\end{document}